\documentclass{article}
\usepackage{amssymb}
\usepackage{latexsym}
\usepackage[latin1]{inputenc}

\usepackage{graphicx}
\usepackage{amsthm}
\usepackage{epsfig}

\input xy
\xyoption{all}
\newdir{ >}{{}*!/-12pt/@{>}}

\newcommand{\eqr}[1]{\rotatebox{#1}{{\tiny $\!\! \sim \!\!\!$}}}

\title{Inductive Lusternik-Schnirelmann category in a model category}
\author{Garc\'{\i}a-Calcines, J.M.* \\[0.3pc]
        Garc\'{\i}a-D\'{\i}az, P.R. }

\date{}

\begin{document}

\maketitle

\footnote{*email address: jmgarcal@ull.es \hfill\break This work
is supported by a grant from MEC (Spain) project MTM2006-06317.
\hfill\break 2000 {\em Mathematics Subject Classification}:
55U35,55M30.
\newline {\em Keywords and phrases}: Lusternik-Schnirelmann category, inductive category,
model category.\hfill\break }

\setlength{\baselineskip}{1.2pc}         
\setlength{\parskip} {0.3pc}         
\newtheorem{theorem}{Theorem}[section]

\newtheorem{lemma}[theorem]{Lemma}
\newtheorem{definition}[theorem]{Definition}
\newtheorem{proposition}[theorem]{Proposition}
\newtheorem{corollary}[theorem]{Corollary}

\newtheorem{remark}[theorem]{Remark}
\newcommand{\Dt}{{\sf Proof:}}
\newcommand{\cqd}{\vspace{-15pt} \hfill$\blacksquare$}

\begin{abstract}
We introduce the notion of inductive category in a model category
and prove that it agrees with the Ganea approach given by
Doeraene. This notion also coincides with the topological one when
we consider the category of (well-) pointed topological spaces.
\end{abstract}

\section*{Introduction}
The Lusternik-Schnirelmann category $\mbox{cat}\hspace{1pt}(X)$ of
a space $X,$ LS-category for short, is a homotopy numerical
invariant which was introduced by the quoted authors in the early
30's in their research on calculus of variations \cite{L-S}. It
has turn out to be an important invariant not only in algebraic
topology but also in other important subjects in mathematics such
as differential topology or dynamical systems. For an excellent
introduction on LS-category theory we refer the reader to
\cite{C-L-O-T} and \cite{J}. Unfortunately, although its
definition is quite simple to establish, the LS-category of a
space is hard to compute. Therefore since its beginnings there
have been different attempts of giving alternative descriptions,
approaches or reasonable bounds in a more algebraic way. There are
four standard formulations of LS-category, at least equivalent for
a large class of topological spaces:

\begin{enumerate}
\item The definition by coverings \cite{F}: that it, the category
$\mbox{cat}\hspace{1pt}(X)$, of a space $X$ is the least n (or
infinite) for which there is a covering of $X$ by n+1 open
subspaces, each of them contractible in $X.$

\item The Whitehead characterization \cite{W}:
$\mbox{cat}\hspace{1pt}(X)\leq n$ if and only if the n+1 diagonal
map can be factorized, up to homotopy, through the n-th fat wedge.

\item
The Ganea characterization \cite{G}:
$\mbox{cat}\hspace{1pt}(X)\leq n$ if and only if the n-th Ganea
fibration (fibre-cofibre construction) admits a homotopy section.

\item The inductive category \cite{G}, \cite{G2}:
$\mbox{indcat}\hspace{1pt}(X)=0$ if and only if $X$ is
contractible; and $\mbox{indcat}\hspace{1pt}(X)\leq n$ if and only
if there exists a cofibration $A\rightarrow Y\rightarrow C$ such
that $\mbox{indcat}\hspace{1pt}(Y)\leq n-1$ and $C$ dominates $X.$
\end{enumerate}

The three latter notions are functorial and also the most usual
and successful ones. These approaches have played an important
role during the development of this invariant and we may assert
that without them it would have been much more difficult to
achieve many important accomplishments given in this subject

At the same time, many direct -upper or lower- bounds of the
LS-category have appeared such as the strong category, the weak
category or the sigma-category. Yet, there is also an important
technique for obtaining lower bounds. It consists of taking models
for topological spaces in some algebraic category where a
LS-category-type homotopy invariant is defined. Necessarily, this
invariant must be established in an abstract homotopy setting
-hopefully a model category in the sense of Quillen \cite{Q},
\cite{Q2}. Then, the algebraic LS-category of the model of X is a
lower bound of the original LS-category of such space. We can
mention the notorious work of Félix and Halperin \cite{F-H}  in
rational homotopy, where they defined a numerical homotopy
invariant, $\mbox{cat} _0,$ in the category of commutative cochain
algebras over $\mathbb{Q}$ and proved that the LS-category of the
rationalization of a space agrees with the $\mbox{cat}_0$ of its
Sullivan model. Halperin and Lemaire \cite{H-L}  also defined
important similar numerical invariants in certain full
subcategories. Thus, different algebraic LS-category-type homotopy
invariants have been appearing in several categories others than
topological spaces.

In order to give a unified theory, basically a generalization of
LS-category in all these categories, Doeraene \cite{D}
successfully introduced an intrinsic notion of LS-category of an
object in a Quillen's model category. Actually, he established
what he called J-category, where cofibrations, fibrations and weak
equivalences take part.  This structure is determined by a certain
set of axioms, which are sufficient to develop an abstract
LS-category theory. Inspired by the topological case he gave two
'a priori' different notions, analogous to the Ganea and Whitehead
characterizations. A crucial point,  the cube axiom, gives the
expected equivalence between these notions. Later, a third
equivalent notion, inspired by the original topological
LS-category notion by coverings, was established by Hess and
Lemaire \cite{He-L}. However, since now, it seems that nobody has
noticed the lack of a notion of abstract inductive category in
this framework, equivalent with the latter.

Our aim in this paper is to give this fourth equivalent notion in
a model category, analogous to that of the topological inductive
category.  In fact, we prove that our abstract inductive category
agrees with the Ganea approach in a J-category.
It is also important to remark that in the case of
topological spaces it coincides with the usual inductive category.

We have divided this article in two sections. In the first one we
give some background about J-categories, necessary for the rest of
the paper. Then, in section 2 we introduce the main notion of this
paper, $\mbox{indcat},$ the inductive category, as well as the
corresponding establishment $\mbox{indcat}\equiv \mbox{cat}.$ For
that goal we firstly set and give some properties of a certain
notion of domination, a bit weaker than having a weak section, as
introduced by Doeraene. Finally we display some examples where our
theory may be of interest.

\section{Preliminaries}

This section is devoted to recall all the notions and results that
will be used along this paper. We begin by giving the definition
of a $J$-category.

A $J$-category is category $\mathcal{C}$ together with a zero
object $0$ and three classes of morphisms called fibrations
($\twoheadrightarrow $), cofibrations ($\rightarrowtail $) and
weak equivalences ($\stackrel{\eqr{0}\:\,}{\rightarrow}$),
satisfying the following set of axioms (J1)-(J5).  Before stating
such axioms some points should be clarified: A morphism which is
both a fibration (resp. cofibration) and a weak equivalence is
called {\it trivial fibration} (resp. {\it trivial cofibration}).
An object $B$ is called {\it cofibrant model} (resp. {\it fibrant
model}) if every trivial fibration
$p:E\stackrel{\eqr{0}\,}{\twoheadrightarrow }B$ admits a section
(resp. if every trivial cofibration
$i:B\stackrel{\eqr{0}\,}{\rightarrowtail }E$ admits a retraction).

\subsection{The axioms of a $J$-category.}

\begin{enumerate}
\item[(J1)] Isomorphisms are trivial
cofibrations and trivial fibrations. The composite of fibrations
(resp. cofibrations) is a fibration (resp. a cofibration). Given
$f:X\rightarrow Y$ and $g:Y\rightarrow Z$ morphisms, if any two of
$f,$ $g,$ $gf$ are weak equivalences then so is the third.

\item[(J2)] For any fibration $p:E\twoheadrightarrow B$ and
morphism $f:B'\rightarrow B$ the pull-back exists in $\mathcal{C}$
$$\xymatrix{
  {E'} \ar@{>>}[d]^{\overline{p}} \ar[r]^{\overline{f}}
                      & E \ar@{>>}[d]^{p}    \\
  {B'} \ar[r]_{f}     & B               }
$$
\noindent and the base extension $\overline{p}$ is a fibration.
Moreover, if $f$ is a weak equivalence then so is $\overline{f};$
and if $p$ is a weak equivalence then so is $\overline{p}.$

Dually, for any cofibration $i:A\rightarrowtail B$ and morphism
$f:A\rightarrow A'$ the push-out of $i$ and $f$ exists, and the
cobase extension $\overline{i}$ of $i$ is a cofibration. If $f$ is
a weak equivalence then so is its cobase extension $\overline{f};$
and if $i$ is a weak equivalence then so is $\overline{i}.$

\item[(J3)] For any map $f:X\rightarrow Y$ there exist
\begin{enumerate}
\item [(i)] an $F$-factorization, that is, a
factorization $f=p\tau $ where $\tau $ is a weak equivalence and
$p$ is a fibration; and

\item[(ii)] a $C$-factorization, that is, a
factorization $f=\sigma i,$ where $i$ is a cofibration and $\sigma
$ is a weak equivalence.
\end{enumerate}

\item[(J4)] Given any object $X$ in $\mathcal{C},$ there exists a
trivial fibration $F\stackrel{\eqr{0}\,}{\twoheadrightarrow }X,$
where $F$ is a cofibrant model.

\medskip For the next axiom, we need the definitions of homotopy
pull-back and homotopy push-out. A commutative square
$$\xymatrix{
  {D} \ar[d]^{g'} \ar[r]^{f'}
                      & C \ar[d]^{g}    \\
  {A} \ar[r]_{f}     & B               }
$$ \noindent is said to be a homotopy pull-back if for some
(equivalently any)
$F$-factorization of $g,$  the induced map 
from $D$ to the pull-back $E'=A\times _B E$ is a weak equivalence
$$\xymatrix@C=.7cm@R=.5cm{
{D} \ar[dd]_{g'} \ar[rrr]^{f'} \ar@{.>}[dr] & & & {C} \ar[dd]^g
\ar[dl]^{\tau }_{\eqr{45}} \\
 & {E'} \ar@{>>}[dl]^{\overline{p}} \ar[r]_{\overline{f}} & {E} \ar@{>>}[dr]^p & \\
 {A} \ar[rrr]_f & & & B  }$$
We can also use an $F$-factorization of $f$ instead of $g$ (or
both). The Eckmann-Hilton dual notion, taking a $C$-factorization
and a push-out, is called homotopy push-out.

\item[(J5)] {\it The cube axiom}. Given any commutative cube
where the bottom face is a homotopy push-out and the vertical
faces are homotopy pull-backs then the top face is a homotopy
push-out.
\end{enumerate}

The fundamental construction which can be made in a $J$-category
is that of the join of two objects over a third.

\begin{definition}
Given two morphisms $f:A\rightarrow B$ and $g:C\rightarrow B$ with
the same target, we consider their join $A*_B C$ as follows. First
consider any $F$-factorization of $g=p\tau $ and the pull-back of
$f$ and $p.$ Let $\overline{f}$ and $\overline{p}$ the base
extensions of $f$ and $p$ respectively. Then take any
$C$-factorization of $\overline{f}=\sigma i$ and the push-out of
$\overline{p}$ and $i.$ This push-out object is denoted by $A*_B
C$ and called the join of $A$ and $C$ over $B.$ The dotted induced
map from $A*_B C$ to $B$ is called the join morphism of $f$ and
$g.$
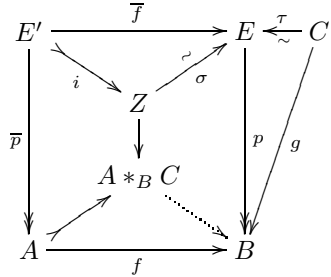
$$\xymatrix@R=.5cm@C=.5cm{
{E'} \ar[rr]^{\overline{f}} \ar@{ >->}[dr]_i
\ar@{>>}[ddd]_{\overline{p}} & & {E} \ar@{>>}[ddd]^p & {C}
\ar[l]_{\tau }^{\eqr{0}} \ar[dddl]^g \\
& {Z} \ar[ur]_{\sigma }^{\eqr{40} } \ar[d] & & \\
& {A*_B C} \ar@{.>}[dr] & & \\
{A} \ar@{ >->}[ur] \ar[rr]_f & & B & }$$
\end{definition}

Two objects $X$ and $Y$ in $\mathcal{C}$ are said to be {\it
weakly equivalent} if there exists a finite chain of weak
equivalences joining $X$ and $Y$
$$\xymatrix{
X \ar@{-}[r]^{\eqr{0}} & \bullet \ar@{-}[r]^(.33){\eqr{0}} &
\bullet \: \cdots \cdots \: {\bullet }\ar@{-}[r]^(.67){\eqr{0}} &
Y }$$ where the symbol $\xymatrix{{\bullet } \ar@{-}[r] &
{\bullet}}$ means an arrow with either left or right orientation.

The object $A*_B C$ and the join map are well defined and
symmetrical up to weak equivalence.

\begin{definition}

Let $f:A\rightarrow B$ and $g:C\rightarrow B$ be morphisms in
$\mathcal{C}.$ We say that $f$ admits a weak lifting along $g$ if
for some (equivalently for any) $F$-factorization $g=p\tau $ there
exists a commutative diagram
$$\xymatrix@R=.4cm@C=.4cm{
 & & C \ar[dd]^g \ar[dl]^{\tau }_{\eqr{45}} \\
 & E \ar@{>>}[dr]^p & \\
 A \ar@{.>}[ur]^s \ar[rr]_f &  & B } $$

In the particular case $f=id_B$ we say that $g:C\rightarrow B$
admits a weak section.
\end{definition}

Now we are ready to the definition of $n$th-Ganea map as well as
the category of a given object $B$ in $\mathcal{C}$.

\begin{definition}
Let $B$ an object in $\mathcal{C}.$ The $n$th Ganea object $G_nB,$
as well as the $n$th Ganea map $p_n^B:G_nB\rightarrow B$ ($n\geq
0$) are constructed as follows: For $n=0$ let $p_0^B$ be the zero
map $0\rightarrow B;$ if $p_{n-1}^B$ is already constructed then
$p_n^B$ is defined as the join map of $0\rightarrow B$ and
$p_{n-1}^B:G_{n-1}\rightarrow B$ (so $G_nB=0*_B G_{n-1}B$).
\end{definition}
\begin{remark}
This construction is functorial. Given $f:B\rightarrow B'$
any morphism in $\mathcal{C}$ we can construct a morphism
$G_n(f):G_nB\rightarrow G_nB'$ such that $p_n^{B'}G_n(f)=fp_n^B.$
\end{remark}

\begin{definition}\cite[3.8]{D}
We say that $\mbox{cat}\hspace{1pt}(B)\leq n$ if and only if the
$n$th Ganea map $p_n:G_nB\rightarrow B$ admits a weak section. If
no such $n$ exists then $\mbox{cat}\hspace{1pt}(B)=\infty .$
\end{definition}
Doeraene proved that this construction is invariant by weak
equivalence, that is, if $B$ and $B'$ are weakly equivalent
objects then
$\mbox{cat}\hspace{1pt}(B)=\mbox{cat}\hspace{1pt}(B').$

We note that, actually, this definition holds in any pointed
category (that is, with zero object) verifying (J1)-(J4) axioms,
without the cube axiom (J5). The cube axiom is just needed to
prove that a second definition, the Whitehead approach, is
equivalent to this one \cite{D}. Since we are just dealing with
the Ganea approach, (J5) axiom will no be needed in this article.
On the other hand, instead of (J3) and (J4) axioms we are
interested in the following slightly stronger ones, given in any
model category.

\begin{enumerate}
\item[(M1)] Given any commutative diagram of solid arrows
$$\xymatrix{
{A} \ar[r] \ar@{ >->}[d]_i & {E} \ar@{>>}[d]^p \\
{X} \ar@{.>}[ur] \ar[r] & {B} }$$ \noindent where $i$ is a
cofibration, $p$ is a fibration and either $i$ or $p$ is a weak
equivalence, then the dotted arrow exists making commutative the
two triangles.

\item[(M2)]
Any map $f:X\rightarrow Y,$ can be factored in two ways:
\begin{enumerate}
\item[(i)] $f=p\tau ,$ where $\tau $ is a trivial cofibration and
$p$ is a fibration ($F$-factorization).

\item[(ii)] $f=\sigma i,$ where $i$ is a cofibration and $\sigma $
is a trivial fibration ($C$-factorization).
\end{enumerate}
\end{enumerate}

From (M2) axiom, (J3) is obviously satisfied. Yet (J4) is also
fulfilled; indeed for any object $X$ we can consider the following
factorization
$$\xymatrix@C=.4cm@R=.4cm{
0 \ar@{ >->}[dr] \ar[rr] & & X \\
 & QX \ar@{>>}[ur]_{p_X}^(.43){\eqr{45} }  }$$
\noindent obtaining, taking into account (M1) axiom, a cofibrant
model.

Dually, we can also consider its fibrant model, obtained by the
factorization of the zero map $X\rightarrow 0$ through a trivial
cofibration $i_X:X\stackrel{\eqr{0}\:}{\rightarrowtail}RX $
followed by a fibration. We also note that given any map
$f:X\rightarrow Y$ it is possible to find a map $Qf:QX\rightarrow
QY$ (respectively $Rf:RX\rightarrow RY$) such that $fp_X=p_YQf$
(respectively $Rfi_X=i_Yf$).

At the sight of these remarks, the framework in which we will be
immersed throughout this paper could be a pointed proper model
category. Nevertheless, the reader may also think that we are
considering a pointed category $\mathcal{C}$ such that (J1),(J2),
(M1) and (M2) axioms are satisfied.

\section{The inductive category.}


We begin by giving the notion of {\it domination}, which is weaker
of that of having a weak section.
\begin{definition} \label{def_dominancia}
Given $X,Y$ objects in $\mathcal{C}$ we will say that $X$ {\it
dominates} $Y$ (denoted $X \gg Y$) if for some (equivalent any)
cofibrant model $QX$ of $X$ and for some (equivalent any) fibrant
model $RY$ of $Y$ there exists a morphism $\alpha:QX\rightarrow
RY$ such that $i_Y:Y\stackrel{\eqr{0}\: }{\rightarrowtail }RY$
admits a weak lifting along $\alpha$
$$\xymatrix@R=.4cm@C=.4cm{
   & & QX \ar[dd]^{\alpha} \ar[dl]_{\eqr{45}}^{\tau} \\
   & E \ar@{>>}[dr]^p &  \\
   Y \ar@{.>}[ur]^s \ar@{ >->}[rr]_{i_Y}^{\sim } &  & RY   }$$
\end{definition}

Note that this definition agrees with the one in the case of
topological spaces. The following properties will be useful for
our notion of inductive category.

\begin{proposition} \label{prop_previo}
${}$
\begin{enumerate}
    \item [i)]  $X\gg X,$ for any object $X$ in $\mathcal{C}.$

    \item [ii)] If $X$ and $Y$ are weakly equivalent then
    $X\gg Z$ if and only if $Y \gg Z.$

    \item [iii)] If $X\stackrel{f}{\rightarrow }Y$ admits a weak section then
    $X\gg Y$.

    \item [iv)] If $X \gg Y$ then $\mbox{cat}\hspace{3pt}{Y}
    \leq \mbox{cat}\hspace{3pt}{X}$.

\end{enumerate}
\end{proposition}

\Dt

\textit{i)} This item is easily proved taking the composite
 $\xymatrix{ QX \ar@{>>}[r]^{\eqr{0} } & X
 \ar@{ >->}[r]^(.45){\eqr{0} } & RX }.$

\textit{ii)} We can suppose, without losing generality, that there
is a weak equivalence $w:X \stackrel{\eqr{0}\:}{\rightarrow} Y.$
Assume that $X\gg Z;$ then there is a morphism $\alpha
:QX\rightarrow RZ$ such that $i_Z$ admits a weak lifting along
$\alpha .$

We can take a cofibrant model of $Y$ in such a way
$Qw:QX\stackrel{\eqr{0}\:}{\rightarrowtail }QY$ is a trivial
cofibration. Since $RZ$ is a fibrant object we can consider a
morphism $\lambda :QY\rightarrow RZ$ such that $\lambda Qw=\alpha
.$ Taking into account that the existence of the weak lifting
along $\alpha $ is independent of the chosen factorization of
$\alpha ,$ we have that $Y\gg Z.$ We have just to take into
account the following diagram:
$$\xymatrix@R=.4cm@C=.65cm{
    & & QY \ar[dd]^{\lambda}
    \ar[dl]_(.52){\eqr{45}}^{\tau} & QX \ar@{.>}[ddl]^{\alpha}
    \ar@{ >->}[l]_(.5){\eqr{0}}^(.5){Qw} \\
    & E' \ar@{>>}[dr]^{p} &  \\
    Z \ar@{.>}[ur] \ar@{ >->}[rr]^{\eqr{0}}_{i_Z} &  & RZ   }$$
The converse is straightforward and left to the reader.

 \textit{iii)} Consider a weak section for $f$
   $$\xymatrix@R=.4cm@C=.5cm{
   & & X \ar[dd]^{f} \ar[dl]_{\eqr{45}}^l \\
   & E \ar@{>>}[dr]_{q} &  \\
   Y \ar@{.>}[ur]^{s} \ar[rr]_{id} &  & Y   }$$
a cofibrant model $QX$ for $X$ and a fibrant model $RY$ for $Y.$
Consider also $\xymatrix{E \ar[r]_(.42){\eqr{0}}^(.42)h & E'
\ar@{>>}[r]^(.42)g & RY}$ an $F$-factorization of $i_Y q.$ Then
the result follows from the following commutative diagram
$$\xymatrix@R=.4cm@C=.4cm{
    & & QX \ar[dd]^{i_Y fp_X} \ar[dl]^{\eqr{45}}_{hlp_X} \\
    & E' \ar@{>>}[dr]^{g} &  \\
    Y \ar@{.>}[ur]^{hs} \ar[rr]_{i_Y} &  & RY   }$$

\textit{iv)} Suppose that $\mbox{cat}\hspace{1pt}{X}=n$ and denote
by $p:E\twoheadrightarrow RY$ and $s:Y \rightarrow E$ as given in
definition \ref{def_dominancia}. Since
$\mbox{cat}\hspace{3pt}{X}=\mbox{cat}\hspace{3pt}{QX}=
\mbox{cat}\hspace{3pt}E$, there exists $\sigma$ a section of
$p_n^E:G_nE\twoheadrightarrow E$. Then we have a lifting
$$\xymatrix@C=.7cm@R=.7cm{
  Y \ar@{ >->}[d]^{\eqr{90}}_{i_Y} \ar[rr]^(.4){G_{n}(p) \sigma s} &  & G_n(RY)
  \ar@{>>}[d]^{\:p_n^{RY}} \\
  RY \ar[rr]_{id} \ar@{.>}[rru] &  & RY   }$$
\noindent showing that
$\mbox{cat}\hspace{1pt}Y=\mbox{cat}\hspace{1pt}RY\leq n.$ \cqd

\bigskip
Now we recall the notion of weak push-out, given in \cite{D}.

\begin{definition}\cite[2.3]{D}
Let $f:A\rightarrow B,$ $f':A'\rightarrow B'$ and $a:A'\rightarrow
A$ be morphisms in $\mathcal{C}.$ We say that
$A'\mbox{-}A\mbox{-}B\mbox{-}B'$ forms a {\it weak push-out} if
for some (equivalently, any) $C$-factorization $\xymatrix{A' \ar@{
>->}[r]^i & X \ar[r]^(.45){\sigma }_(.45){\eqr{0} } & B' }$ of $f',$ one has a
homotopy push-out $A'\mbox{-}A\mbox{-}B\mbox{-}X$
$$\xymatrix{
A' \ar[rr]^a \ar@{ >->}[dr]^i \ar[d]_{f'} & & A \ar[d]^f \\
B' & X \ar[l]^(.46){\sigma }_(.46){\eqr{0}} \ar[r]_x  & B}$$

We say that $f$ is the {\it weak cobase extension} of $f'$ by $a.$
\end{definition}
In particular, if $a:A'\rightarrow 0$ is the zero map we say that
the map $x$ (or just the object $B$) is the {\it homotopy cofibre}
of $f':A'\rightarrow B'.$

Given $f':A'\rightarrow B'$ and $a:A'\rightarrow A$ morphisms in
$\mathcal{C},$ the usual way to obtain their weak push-out (up to
weak equivalence) is to consider a $C$-factorization $f'=\sigma i$
and then the push-out of $i$ along $a.$ Considering the gluing
lemma \cite{B}, this construction is well-defined and symmetrical
up to weak equivalence (i.e. we may take a $C$-factorization
$a=\sigma i$ and then form the push-out of $i$ along $f'.$)

Then, given $f:A\rightarrow Y$ morphism in $\mathcal{C},$ its {\it
cofibre sequence} $A\stackrel{f}{\rightarrow
}Y\stackrel{p}{\rightarrow} C$ might be obtained as the following
push-out
$$\xymatrix{
 & A \ar[r]^f \ar@{ >->}[d]^k \ar[dl] & Y
 \ar@{ >->}[d]^{\overline{k}=p} \\
 0 & CA \ar[l]^{\eqr{0}} \ar[r]_{\overline{f}} & {C} }$$

Now we are giving the main notion of this paper.
\begin{definition}
Let $X$ be any object in $\mathcal{C}.$ The inductive category of
$X$, $\mbox{indcat}\hspace{3pt}{X}$, is defined as follows:
$\mbox{indcat}\hspace{3pt}{X}=0$ if and only if the zero map
$0\longrightarrow X$ admits a weak section; for $n\geq 1,$
$\mbox{indcat}\hspace{3pt}{X}\leq n$ if and only if there exists a
cofibre sequence
$$A\stackrel{f}{\longrightarrow }Y\stackrel{p}\longrightarrow C$$
\noindent such that $\mbox{indcat}\hspace{2pt}{Y}\leq n-1$ and
$C\gg X.$
\end{definition}

\begin{remark}
Considering proposition \ref{prop_previo} it is straightforward to
check that $\mbox{indcat}\hspace{3pt}X$ is well defined and it is
invariant up to weak equivalence.
\end{remark}

\begin{theorem}
Given $X$ any object in $\mathcal{C}$ we have
$$\mbox{indcat}\hspace{2pt}{X}=\mbox{cat}\hspace{2pt}{X}$$
\end{theorem}

\Dt

By using inductive arguments and proposition \ref{prop_previo} we
can plainly see that $\mbox{indcat}\hspace{2pt}{G_k(X)}\leq k,$
for all $k\geq 0.$  Therefore, and again by the same proposition
\ref{prop_previo}, we have $\mbox{indcat}\hspace{2pt}{X}\leq
\mbox{cat}\hspace{2pt}{X}$.

For the inequality $\mbox{cat}\hspace{2pt}{X}\leq
\mbox{indcat}\hspace{2pt}{X}$ we proceed by induction on the
integer $\mbox{indcat}\hspace{2pt}{X}=k$ and the object $X.$  The
result is obvious for $k=0$. Now suppose that the statement is
true for any object $Z$ and $k\leq n-1$ and that
$\mbox{indcat}\hspace{2pt}{X}=n$. Then, there is a cofibre
sequence $A \stackrel{f}{\longrightarrow} Y
\stackrel{p}{\longrightarrow} C$ such that $C\gg X$ and
$\mbox{cat}\hspace{2pt}{Y}=\mbox{indcat}\hspace{2pt}{Y}\leq n-1.$

We explicitly give a construction, up to weak equivalence, of the
$n$th Ganea map of $C$ as explained in the next diagram:
$$\xymatrix@R=0.4cm@C=0.9cm{
  F_{n-1}(C) \ar@{ >->}[dr]^{\varepsilon}
  \ar[ddd]^{\overline{\overline{f}}} \ar@{>>}[rr]^{\overline{p^{C}_{n-1}}} & &
  CA \ar[ddd]^{\overline{f}} \\
   & \widehat{CA} \ar[d]_{\mu} \ar@{>>}[ur]^(.45){q}_(.45){\eqr{45}} &  \\
   & G_n(C) \ar@{.>}[dr]^{p^{C}_n}  &  \\
   G_{n-1}(C) \ar@{>>}[rr]_{p^{C}_{n-1}} \ar[ur]^{w} &  & {C}   }$$
Observe that we can assume that
$p^C_{n-1}:G_{n-1}C\twoheadrightarrow C$ is already a fibration.
Then the homotopy fibre $F_{n-1}(C)$ of $p^C_{n-1}$ may be
obtained, up to weak equivalence, as the pull-back of $p^C_{n-1}$
and $\overline{f}:CA\rightarrow C.$ Next we factor
$\overline{p^{C}_{n-1}},$ the base change of $p^{C}_{n-1},$ as a
cofibration $\varepsilon $ followed by a trivial fibration $q.$
Finally, take the push-out $G_n(C)$ of $\varepsilon $ and
$\overline{\overline{f}}$ as well as the push-out map $p_n^C.$

Considering a section $s:Y \rightarrow G_{n-1}(Y)$ of
$p_{n-1}^Y:G_{n-1}(Y)\twoheadrightarrow Y$ and the map
$G_{n-1}(p):G_{n-1}(Y) \rightarrow G_{n-1}(C)$ the Ganea
construction of $p$ we can take the pull-back map
$$\lambda=(G_{n-1}(p)sf,k):A \rightarrow F_{n-1}(C).$$

By (M1) axiom take $q'$ any lift in the diagram
$$\xymatrix@C=1.3cm{
  A \ar@{ >->}[d]_{k} \ar[r]^{\varepsilon \lambda} &
  \widehat{CA} \ar@{>>}[d]_{\eqr{90}}^q \\
  CA \ar@{=}[r] \ar@{.>}[ur]^{q'} & CA   }$$
Then one can straightforwardly check that the push-out map
$$\sigma=(\mu q', wG_{n-1}(p) s):C \rightarrow G_n(C)$$ is a
section of $p_n^C.$ Following proposition \ref{prop_previo}.iv) we
have the desired result \cqd

\ \par
\bigskip
\begin{remark}
We point out that all these notions and results have their dual in
the sense of Eckmann-Hilton. In the obvious way the notion of {\it
inductive cocategory} can be defined; and using the dual results
we have that the inductive cocategory agrees with the Ganea
approach of cocategory (by using cojoins).

\end{remark}

\subsection{Some examples}

Now we review some examples of pointed categories in which the
inductive category can be applied. Of course, any pointed proper
model category may be taken.

\begin{itemize}
    \item [(i)] $\mathbf{Top}^w,$ the category of well-pointed
topological spaces and continuous maps which preserve the base
point.

By a well-pointed space we mean a pointed space $X$ in which the
inclusion of the base point in $X$ is a closed topological
cofibration (that is, a closed map such that verifies the homotopy
extension property). Considering
\begin{itemize}
    \item Fibrations: the Hurewicz fibrations,
    \item Cofibrations: the closed cofibrations,
    \item Weak equivalences: the homotopy equivalences,
\end{itemize}
then $\mathbf{Top}^w$ together these classes of maps satisfies
(J1), (J2), (M1) and (M2) axioms and all spaces are fibrant and
cofibrant. In fact, Str\o m proved that it verifies stronger
conditions \cite{S}.

As it is well-known, T. Ganea \cite{G}, \cite{G2} proved that the
inductive (co)category agrees with the Ganea approach of
Lusternik-Schnirelmann (co)category in the context of topological
spaces. Here we have displayed a different and slightly simpler
proof of this fact.


    \item [(ii)] The category $\mathbf{S}^{\bullet }$ of pointed simplicial
sets and simplicial maps preserving base points, where
\begin{itemize}
    \item Fibrations: Kan fibrations,
    \item Cofibrations: injective simplicial maps,
    \item Weak equivalences: maps whose geometric realizations are homotopy
equivalences.
\end{itemize}
Then, $\mathbf{S}^{\bullet }$ is a proper model category \cite{Q}
where all objects are cofibrant and the fibrant objects are the
Kan complexes. We can consider Kan's $\mathbf{Ex}^{\infty }$
functor which associates a Kan complex $\mathbf{Ex}^{\infty
}\hspace{2pt}L$ to any simplicial set $L$ up to weak equivalence.
This construction also involves a simplicial map $\nu
_L:L\stackrel{\eqr{0}\, }{\rightarrowtail }\mathbf{Ex}^{\infty
}\hspace{2pt}L,$ which is, actually, a fibrant model for $L.$
Then, given $K,$ $L$ pointed simplicial sets, $K$ dominates $L$ if
there exists a simplicial map $\alpha :K\rightarrow
\mathbf{Ex}^{\infty }\hspace{2pt}L$ such that $\nu _L$ admits a
weak section along $\alpha $
$$\xymatrix@R=.4cm@C=.4cm{
   & & K \ar[dd]^{\alpha} \ar[dl]_{\eqr{45}}^{\tau} \\
   & K' \ar@{>>}[dr]^p &  \\
   L \ar@{.>}[ur]^s \ar@{ >->}[rr]_(.45){\nu _L}^(.43){\eqr{0} } &  &
   {\mathbf{Ex}^{\infty }\hspace{2pt}L}   }$$

    \item [(iii)] The category $\mathbf{CDA^*}$ of augmented
    conmutative cochain algebras over a field of
    characteristic zero. Considering
    \begin{itemize}
        \item Fibrations: Surjective maps.
        \item Cofibrations: KS-extensions.
        \item Weak equivalences: Quasi-isomorphism, that is, maps which induce isomorphisms in
        cohomology.
    \end{itemize}
    and the full subcategory $\mathbf{CDA^*}^{\mathrm{c0}}$ of
    c-connected differential algebras. We can find in \cite{H-W} and \cite{B}
    that $\mathbf{CDA^*}^{\mathrm{c0}}$ with the induced structure satisfies
    (J1), (J2), (M1) and (M2) axioms, where all objects
    are fibrant. The $n$th Ganea algebra of $\Lambda X,$ a minimal c-1-connected
    KS-complex, defined by Félix and Halperin \cite{F-H}
    was interpreted by Doeraene in another way as
    a cojoin construction \cite{D}. That is, the dual $n$th
    Ganea map $\Lambda X \rightarrow G^n(\Lambda X)$ admits a weak
    retraction if and only if the projection $\Lambda X
    \rightarrow (\Lambda X / \Lambda^{>n} X)$ admits a weak
    retraction. So in this case the Doeraene's notion of cocat in
    $\mathbf{CDA^*}^{\mathrm{c0}}$ agrees with the rational LS-category
    defined by Félix and Halperin. We have another equivalent definition in
    terms of inductive cocategory:

    We say that $\mbox{indcocat}\hspace{3pt}{\Lambda X} \leq n$ if
    there is a fiber sequence (in $\mathbf{CDA^*}^{\mathrm{c0}}$)
    $$F\stackrel{}{\longrightarrow }E\stackrel{}\longrightarrow B$$
    with $\mbox{indcocat}\hspace{3pt}E \leq n-1$
    and there exists a map $\Lambda V \rightarrow \Lambda X$
    which admits a weak section, being $\Lambda V $ a minimal
    model for $F.$
\end{itemize}

\end{document}